\documentclass[11pt]{article}

\usepackage{amsmath,amsfonts,amsthm,amssymb,tikz,tikz-cd,url,hyperref,cleveref,old-arrows,mathtools}

\pdfoutput=1

\theoremstyle{plain}
\newtheorem{thm}{Theorem}[section]

\newtheorem{lem}[thm]{Lemma}
\newtheorem{prop}[thm]{Proposition}

\newtheorem*{thmnonum}{Theorem}

\theoremstyle{definition}

\newtheorem{rem}[thm]{Remark}

\theoremstyle{remark}

\newcommand{\bbB}{\mathbb{B}}
\newcommand{\bbC}{\mathbb{C}}

\newcommand{\bbF}{\mathbb{F}}

\newcommand{\bbH}{\mathbb{H}}

\newcommand{\bbN}{\mathbb{N}}

\newcommand{\bbP}{\mathbb{P}}
\newcommand{\bbQ}{\mathbb{Q}}
\newcommand{\bbR}{\mathbb{R}}

\newcommand{\bbZ}{\mathbb{Z}}

\newcommand{\bfG}{\mathbf{G}}

\newcommand{\calC}{\mathcal{C}}

\newcommand{\calL}{\mathcal{L}}
\newcommand{\calM}{\mathcal{M}}

\newcommand{\calO}{\mathcal{O}}

\newcommand{\frakU}{\mathfrak{U}}

\newcommand{\frakg}{\mathfrak{g}}

\newcommand{\frakp}{\mathfrak{p}}

\newcommand{\al}{\alpha}
\newcommand{\gam}{\gamma}
\newcommand{\Gam}{\Gamma}

\newcommand{\de}{\delta}
\newcommand{\Del}{\Delta}

\newcommand{\lam}{\lambda}
\newcommand{\Lam}{\Lambda}
\newcommand{\sig}{\sigma}

\DeclareMathOperator{\SL}{SL}
\DeclareMathOperator{\PSL}{PSL}
\DeclareMathOperator{\GL}{GL}
\DeclareMathOperator{\PGL}{PGL}

\DeclareMathOperator{\PU}{PU}
\DeclareMathOperator{\SU}{SU}

\DeclareMathOperator{\Ad}{Ad}
\DeclareMathOperator{\Tr}{Tr}

\DeclareMathOperator{\Aut}{Aut}

\newcommand{\bs}{\backslash}

\newcommand{\lra}{\longrightarrow}

\newcommand{\ssm}{\smallsetminus}

\newcommand{\conj}{\overline}
\newcommand{\wh}{\widehat}
\newcommand{\wt}{\widetilde}

\newenvironment{pf}{\begin{proof}}{\end{proof}}

\newenvironment{enum}{\begin{enumerate}}{\end{enumerate}}

\allowdisplaybreaks

\usetikzlibrary{arrows}

\title{Finite groups and complex projective surfaces}
\author{Alexander Lubotzky\footnote{Weizmann Institute of Science, \texttt{alex.lubotzky@mail.huji.ac.il}} \and Matthew Stover\footnote{Temple University, \texttt{mstover@temple.edu}}}
\date{\today}

\begin{document}

\maketitle

\begin{abstract}
In response to a question raised by Belolipetsky and the first author in \cite{BelolipetskyLubotzky}, we prove that for every finite group $G$ there are infinitely many isomorphism classes of compact complex hyperbolic $2$-manifolds with automorphism group isomorphic to $G$.
\end{abstract}

\section{Introduction}\label{sec:Intro}

Throughout this paper, an automorphism of a complex manifold will always mean a biholomorphism. The main result of this paper, proved in \Cref{sec:Main}, is the following.

\begin{thm}\label{thm:Main}
For every finite group $G$, there are infinitely many isomorphism classes of closed complex hyperbolic $2$-manifolds with automorphism group isomorphic to $G$.
\end{thm}

Greenberg \cite{Greenberg} proved \Cref{thm:Main} for Riemann surfaces, i.e., complex hyperbolic $1$-manifolds. Belolipetsky and the first author \cite{BelolipetskyLubotzky} proved the analogue of \Cref{thm:Main} for isometry groups of hyperbolic $n$-manifolds, $n \ge 3$, and this paper answers a question raised in \cite[\S 6.5]{BelolipetskyLubotzky}. Our argument follows the broad outline of that paper, but differences at key steps create many new challenges we must overcome.

Recall that a compact complex hyperbolic $n$-manifold is a quotient $\Gam \bs \bbB^n$ of the unit ball in $\bbC^n$ by a torsion-free cocompact lattice $\Gam$ in the Lie group $\PU(n,1)$. Note that $\bbB^1$ is simply the Poincar\'e disk with the hyperbolic metric. For $n \ge 2$, Mostow rigidity then implies that the automorphism group of a complex hyperbolic manifold is precisely the holomorphic isometry group, and the proof of \Cref{thm:Main} is through this interpretation.

One strong similarity between the proofs is in the input: a cocompact nonarithmetic lattice in the isometry group of the symmetric space. Moreover, \cite{BelolipetskyLubotzky} uses a nonabelian free quotient to begin the construction, whereas we use a quotient that is a cocompact Fuchsian group. The subgroup counting methods in both papers are nearly identical. Note that it is quite easy, given a free or surface group quotient of a finite index subgroup $\Lam$ of a lattice $\Gam$, to produce a finite index subgroup $B$ of $\Gam$ so that $N_\Gam(B) / B$ \emph{contains} $G$, where $N_\Gam(B)$ is the normalizer of $B$ in $\Gam$, but ensuring that $B$ is torsion-free and that $N_\Gam(B) / B$ is \emph{exactly} $G$ is quite subtle.

The nonarithmetic lattices used in \cite{BelolipetskyLubotzky} are the famous examples of Gromov and Piatetski-Shapiro \cite{GPS}. These lattices naturally have the structure of a free product with amalgamation over the fundamental group of a codimension one totally geodesic subspace, and the nonabelian free quotient of a finite index subgroup is derived from this (see \cite{LubotzkyFreeQuo}). In the complex hyperbolic setting, there are no real codimension one totally geodesic subspaces \cite[\S 3.1]{Goldman}, hence no such splittings.

Our replacement, from the theory of K\"ahler groups, is the \emph{universal homomorphism} to a direct product of cocompact Fuchsian groups; see \cite{HLPSV} for an introduction. Cocompact nonarithmetic lattices in $\PU(2,1)$ with nontrivial universal homomorphism are found among the famous examples of Livn\'e \cite{Livne} and Deligne--Mostow \cite{DeligneMostow}; see \Cref{thm:Examples}. With this different input, a key lattice $\Del$ constructed early in the outline of \cite{BelolipetskyLubotzky} now fails to be normal in the maximal lattice $\Gam$ in $\PU(2,1)$. In a certain sense, much of the work in \Cref{sec:Main} is using properties of the universal homomorphism to circumvent this significant obstruction to executing several steps in the strategy from \cite{BelolipetskyLubotzky}.

\begin{rem}\label{rem:Highern}
Whether or not there are cocompact nonarithmetic lattices in $\PU(n,1)$ for all $n \ge 3$ is one of the major open problems in the study of lattices in Lie groups. For any $n \ge 3$, if there is a cocompact nonarithmetic lattice in $\PU(n,1)$ with nontrivial universal homomorphism, then this paper produces infinitely many isomorphism classes of closed complex hyperbolic $n$-manifolds with isometry group any finite group $G$ \emph{verbatim}.
\end{rem}

In this paper, a genus $g$ surface group will always mean the fundamental group of a closed Riemann surface of genus $g$, despite the slight clash of terminology with fundamental groups of complex surfaces also appearing. In \Cref{ssec:FuchsianAut} we will prove the following result on automorphisms of surface groups, which is an analogue of a theorem of the first author in the free group case \cite[Thm.\ 1]{LubotzkyFree} that may be of independent interest. 

\begin{thm}\label{thm:Fixpsubs}
Let $T$ be a surface group of genus $g \ge 2$, $p$ be a prime, and suppose that $\varphi \in \Aut(T)$ has the property that $\varphi(N) = N$ for all finite index normal subgroups of $T$ with index a power of $p$. Then $\varphi$ is inner.
\end{thm}

\noindent\textbf{Interpretation as a complex projective surface}

\medskip

Compact complex hyperbolic $2$-manifolds are also minimal smooth complex projective surfaces of general type \cite[\S V.20]{BPV}. The appendix to this paper recasts \Cref{thm:Main} in that context to prove the following result.

\begin{thm}\label{thm:ProjectiveVersion}
For every finite group $G$, there are infinitely many isomorphism classes of irreducible, minimal smooth complex projective surfaces of general type with automorphism group $G$.
\end{thm}

A surface is irreducible if it is not finitely (\'etale) covered by a product of complex curves. Greenberg's theorem easily produces products of complex curves with any finite automorphism group, so only the irreducible case is of interest. As described in the appendix, one may be able to construct irreducible examples using other methods. For example, there may exist a nonsingular fibration $X \to C$ over a curve $C$ (e.g., a ruled surface) where $\Aut(C) \cong G$ by Greenberg and there is a natural isomorphism $\Aut(X) \cong \Aut(C)$. However, we were not able to find a reference for \Cref{thm:ProjectiveVersion} in the literature.

\begin{rem}\label{rem:HardPart}
Our construction is also related to a holomorphic map to a smooth projective curve (see \Cref{ssec:KahlerGroups}). However, the bulk of the proof is precisely in showing that our surface has no extra automorphisms, which is highly unclear from the nature of the construction.
\end{rem}

\medskip

\noindent\textbf{Organization of the paper}

\medskip

\noindent
This paper is organized as follows. \Cref{sec:Prelim} contains preliminaries about lattices in $\PU(2,1)$, universal homomorphisms for K\"ahler groups, and the proof of \Cref{thm:Fixpsubs}. \Cref{sec:Main} proves \Cref{thm:Main}. The appendix discussions automorphism groups of smooth complex projective surfaces.

\medskip

\noindent \textbf{Acknowledgments}

\medskip

\noindent
Much of the work on this paper was done while the authors were at the conference `Discrete and profinite groups' at the Isaac Newton Institute for Mathematical Sciences, and thank the institute for its hospitality. The first author was supported by the European Research Council (ERC) under the European Union’s Horizon 2020 (N.\ 882751). The second author was supported by Grants DMS-2203555 and DMS-2506896 from the National Science Foundation and the Simons Foundation [SFI-MPS-TSM-00014184, MS].

\section{Preliminaries}\label{sec:Prelim}

This section reviews preliminaries needed for the sections that follow.

\subsection{Lattices in $\PU(2,1)$}\label{ssec:PU21}

Let $\bbB^2$ denote the unit ball in $\bbC^2$ with its Bergman metric of constant holomorphic sectional curvature $-1$. Then the group of holomorphic isometries of $\bbB^2$ can be identified with the Lie group $\PU(2,1)$. If ${\Lam < \PU(2,1)}$ is a discrete, cocompact, and torsion-free subgroup, then $\Lam \bs \bbB^2$ is a complex manifold, and in this case the fact that $\Lam \bs \bbB^2$ is a smooth projective variety goes back to Kodaira \cite{Kodaira}.

By the famous Margulis dichotomy, a lattice $\Lam$ in $\PU(2,1)$ is \emph{arithmetic} if and only if its commensurator
\begin{equation}\label{eq:Commensurator}
\calC(\Lam) \coloneqq \left\{g \in \PU(2,1)\ :\ \Lam \textrm{ and } g \Lam g^{-1} \textrm{ are commensurable}\right\}
\end{equation}
is topologically dense in $\PU(2,1)$ \cite[p.\ 4]{Margulis}. The existence of cocompact nonarithmetic lattices was first established by Mostow \cite{Mostow}, and it is well-known that if $\Lam$ is nonarithmetic then $\calC(\Lam)$ is the \emph{unique} maximal lattice in $\PU(2,1)$ commensurable with $\Lam$.

Every lattice $\Lam < \PU(2,1)$ is \emph{integral} in the following sense. Let $k$ be the \emph{adjoint trace field}, i.e., the field over $\bbQ$ generated by $\{\Tr(\Ad \lam) : \lam \in \Lam\}$, where $\Ad$ is the adjoint representation on the Lie algebra $\frakg$ of $\PU(2,1)$. Let $\calO_k$ denote the ring of integers of $k$. Then there is a vector space $V$ over $k$, a $\calO_k$-lattice $\calL \subset V$ (i.e., an $\calO_k$-submodule of full rank), and an isomorphism $V \otimes_k \bbR \cong \frakg$ for some real embedding of $k$ so that
\begin{equation}\label{eq:Integral}
\Ad(\Lam) \le \Aut(\calL) \le \GL(\frakg).
\end{equation}
In particular, there is a linear algebraic group $\bfG$ defined over $k$ and a structure on $\bfG$ as a group scheme over $\calO_k$ (defined via the action on $\calL$) so that $\Lam < \bfG(\calO_k)$. See \cite[Thm.\ 1.5]{BFMS2} and \cite[Thm.\ 1 \& 3]{Vinberg} for details.

Let $\wt{\bfG}$ be the simply connected algebraic group over $k$ associated with the preimage $\wt{\Lam}$ of $\Lam$ in the simply connected real algebraic group $\SU(2,1)$. Given a prime ideal $\frakp$ of $\calO_k$ with residue field $\bbF_\frakp$, the restriction of reduction
\begin{equation}\label{eq:ReductionModp}
\pi_\frakp : \wt{\bfG}(\calO_k) \lra \wt{\bfG}(\bbF_\frakp)
\end{equation}
to any $\wt{\Del} \le \wt{\bfG}(\calO_k)$ has kernel the \emph{congruence subgroup} $\wt{\Del}(\frakp) \trianglelefteq \wt{\Del}$. Strong approximation for Zariski dense subgroups of simply connected algebraic groups says that if $\wt{\Del}$ is Zariski dense in $\wt{\bfG}(k)$, then $\pi_\frakp(\wt{\Del}) = \wt{\bfG}(\bbF_\frakp)$ for all but finitely many $\frakp$ \cite[Thm.\ 1.1]{Weisfeiler}; also see Cor.\ 3 in Window 9 of \cite{LubotzkySegal} along with discussion on p.\ 392 implying the more general statement for number fields other than $\bbQ$. If $\wt{\Del}$ is the preimage of $\Del$ in $\SU(2,1)$, then congruence subgroups of $\Del$ are those whose preimage in $\wt{\Del}$ are congruence subgroups.

Since $\PU(2,1)$ has type ${}^2 \mathrm{A}_2$, there is a quadratic extension $\ell$ of $k$ canonically associated with $\bfG$. If $\Del \le \Gam$ is Zariski dense in $\PU(2,1)$, then away from a finite set of exceptional primes, the group $\Del / \Del(\frakp)$ has commutator subgroup $\PSL_3(\bbF_\frakp)$ if $\frakp$ splits in $\ell$ or $\PU(3; \bbF_\frakp)$ if $\frakp$ is inert in $\ell$. In particular, the congruence quotient of $\Del$ has commutator subgroup a nonabelian finite simple group for all but finitely many primes $\frakp$. For example in the split case, observe that $\Del / \Del(\frakp)$ is then naturally a subgroup of $\PGL_3(\bbF_\frakp)$ containing $\PSL_3(\bbF_\frakp)$, so it has trivial center.

\subsection{K\"ahler groups}\label{ssec:KahlerGroups}

A \emph{K\"ahler group} is the fundamental group of a compact K\"ahler manifold. See \cite{KahlerBook, Py} for basics on K\"ahler groups. Note that finite index subgroups of K\"ahler groups are naturally K\"ahler groups.

Let $\Gam < \PU(2,1)$ be a cocompact lattice, possibly with torsion. Since $\PU(2,1)$ is a linear Lie group and $\Gam$ is finitely generated, $\Gam$ is residually finite and thus there is a finite index subgroup $\Lam < \Gam$ that is torsion-free. Since smooth projective varieties are compact K\"ahler manifolds, by pullback of the Fubini--Study metric on $\bbP^N$ under a projective embedding, it follows that any such $\Lam$ is a K\"ahler group. It then follows from \cite[Lem.\ 1.1.15]{KahlerBook} that all cocompact lattices in $\PU(2,1)$ are K\"ahler groups.

\medskip

Fundamental to this paper is the \emph{universal homomorphism} attached to any K\"ahler group. See \cite[\S 2.C]{HLPSV} for discussion and references for the following result defining the universal homomorphism.

\begin{thm}\label{thm:Universal}
Suppose $\Lam$ is a K\"ahler group. Then there is a (possibly empty) finite collection of cocompact Fuchsian groups $\{T_\al\}_{\al \in \frakU}$ and a homomorphism
\begin{equation}\label{eq:UniversalHom}
\rho : \Lam \lra T \coloneqq \prod_{\al \in \frakU} T_\al
\end{equation}
such that:
\begin{enum}

\item For each $\al \in \frakU$ with projection $p_\al : T \to T_\al$, the composition $p_\al \circ \rho$ is surjective.

\item Every surjective homomorphism from $\Lam$ to a cocompact Fuchsian group factors through some $p_\al \circ \rho$.

\item A surjective homomorphism $\sig : \Lam \to T^\prime$ onto a cocompact Fuchsian group $T^\prime$ has finitely generated kernel if and only if there is some $\al \in \frakU$ and an isomorphism $\psi : T_\al \to T^\prime$ so that $\sig = \psi \circ p_\al \circ \rho$.

\end{enum}
\end{thm}

The $T_\al$ are called the \emph{direct factors} of $\rho$. When $\Lam$ has no homomorphisms onto cocompact Fuchsian groups, note that the universal homomorphism is simply the trivial homomorphism to the trivial group. When $\frakU \neq \varnothing$, we will say that $\Lam$ has \emph{nontrivial} universal homomorphism. Moreover, we further partition the factors of the universal homomorphism as $\frakU = \frakU_{\mathrm{tf}} \sqcup \frakU_{\mathrm{orb}}$, where $\frakU_{\mathrm{tf}}$ is the (possibly empty) set of $\al \in \frakU$ so that $T_\al$ is a torsion-free Fuchsian group, i.e., a genus $g$ surface group for some $g \ge 2$. Then define
\begin{equation}\label{eq:UniversalTFHom}
\rho_{\mathrm{tf}} : \Lam \lra \prod_{\al \in \frakU_{\mathrm{tf}}} T_\al
\end{equation}
to be the \emph{torsion-free universal homomorphism}. The torsion-free case of the following simple lemma will prove important later.

\begin{lem}\label{lem:UniversalChar}
Let $\Lam$ be a K\"ahler group, $\rho$ be its universal homomorphism, and $\rho_{\mathrm{tf}}$ be the torsion-free universal homomorphism for $\Lam$. Then $\ker(\rho)$ and $\ker(\rho_{\mathrm{tf}})$ are characteristic subgroups of $\Lam$.
\end{lem}
\begin{pf}
If $\frakU$ or $\frakU_{\mathrm{tf}}$ is empty, then the respective conclusion is immediate. Now suppose $\al \in \frakU$ with factor $T_\al$ and $\rho_\al \coloneqq p_\al \circ \rho$ is the homomorphism from $\Lam$ onto $T_\al$. Given any $\varphi \in \Aut(\Lam)$, $\rho_\al \circ \varphi$ is a homomorphism onto $T_\al$ with finitely generated kernel $\varphi^{-1}(\ker(\rho_\al))$. Thus there is a $\varphi(\al) \in \frakU$ and an isomorphism $\psi : T_\al \to T_{\varphi(\al)}$ so that $\rho_\al \circ \varphi = \psi \circ \rho_{\varphi(\al)}$, and therefore
\begin{equation}\label{eq:AutInvariant}
\varphi^{-1}(\ker(\rho_\al)) = \ker(\rho_\al \circ \varphi) = \ker(\rho_{\varphi(\al)}).
\end{equation}
Note that $\varphi(\al) \in \frakU_{\mathrm{tf}}$ if and only if $\al \in \frakU_{\mathrm{tf}}$. Considering $\rho_{\varphi^{-1}(\beta)}$, every $\beta \in \frakU$ is of the form $\varphi(\al)$ for some $\al \in \frakU$. Thus
\begin{align*}
\varphi(\ker(\rho)) &= \bigcap_{\al \in \frakU} \varphi(\ker(\rho_{\varphi(\al)})) \\
&= \bigcap_{\al \in \frakU} \ker(\rho_{\varphi(\al)}) \\
&= \bigcap_{\beta \in \frakU} \ker(\rho_\beta) \\
&= \ker(\rho)
\end{align*}
and similarly for the intersection over $\al \in \frakU_{\mathrm{tf}}$. This proves that $\ker(\rho)$ and $\ker(\rho_{\mathrm{tf}})$ are indeed characteristic in $\Lam$.
\end{pf}

The next simple lemma will also be of use.

\begin{lem}\label{lem:TFFuchsian}
Let $\Lam$ be a K\"ahler group with nontrivial universal homomorphism. Then there is a finite index normal subgroup $\Lam^\prime \le \Lam$ such that the universal homomorphism for $\Lam^\prime$ has a direct factor that is a torsion-free Fuchsian group.
\end{lem}
\begin{pf}
Fix a factor $T_\al$ of the universal homomorphism for $\Lam$ and a homomorphism $\sig : T_\al \to F$ onto a finite group $F$ with $T^\prime_\sig \coloneqq \ker(\sig)$ torsion-free. Then $\Lam^\prime \coloneqq \ker(\sig \circ \rho_\al)$ is a normal finite index subgroup mapping onto the torsion-free Fuchsian group $T^\prime_\sig$ with kernel equal to $\ker(\rho_\al)$, which is finitely generated. Thus $T^\prime_\al$ is isomorphic to a factor of the universal homomorphism for $\Lam^\prime$, as desired.
\end{pf}

Notice that the subgroup $\Lam^\prime$ in \Cref{lem:TFFuchsian} may have additional factors for its universal homomorphism beyond those naturally induced from $\Lam$. Finally, we note that there are cocompact nonarithmetic lattices in $\PU(2,1)$ with nontrivial universal homomorphism. The following is known to experts, but we recount the argument.

\begin{thm}\label{thm:Examples}
There are cocompact nonarithmetic lattices $\Gam$ in $\PU(2,1)$ with nontrivial universal homomorphism.
\end{thm}
\begin{pf}
We use a lattice in $\PU(2,1)$ constructed by Livn\'e \cite{Livne} well-known by experts to have nontrivial universal homomorphism. See \cite[\S 16]{DeligneMostow2} for the following description of his construction. Consider the classical congruence subgroup $\Gam(n)$ of $\SL_2(\bbZ)$ with level $n \in \bbN$ acting on the upper half-plane $\bbH^2$ and the modular curve $X(n)^0 \coloneqq \Gam(n) \bs \bbH^2$. We then have the elliptic modular surface
\begin{equation}\label{eq:Elliptic}
E(n)^0 \coloneqq\!\left(\Gam(n) \ltimes \bbZ^2\right)\! \bs \!\left(\bbH^2 \times \bbC\right)
\end{equation}
where the action of the semidirect product defined by
\begin{equation}\label{eq:EllipticAction}
\left(\begin{pmatrix} a & b \\ c & d \end{pmatrix}\, ,\, (x,y) \right)\!(\tau, z) \coloneqq\! \left( \frac{a \tau + b}{c \tau + d}\, ,\, \frac{z + x \tau + y}{c \tau + d} \right)
\end{equation}
gives a fibration $p : E(n)^0 \to X(n)^0$ with fiber above the point $[\tau] \in X(n)^0$ naturally isomorphic to the elliptic curve $A_\tau \coloneqq \langle 1, \tau \rangle \bs \bbC$.

Now consider $X(n) \coloneqq \Gam(n) \bs (\bbH^2 \times (\bbQ \cup \{\infty\}))$ obtained from one-point compactification at each cusp of the modular curve. There is a compactification $E(n)$ of $E(n)^0$ obtained by adding an $n$-cycle of complex projective lines over each cusp of $X(n)^0$. Specifically, $p$ extends to a natural surjective holomorphic map $p : E(n) \to X(n)$ where a cusp point $z_\infty \in X(n)$ has preimage a union of $n$ copies of the complex projective line $\bbP^1$ whose intersection graph is an $n$-gon. See \Cref{fig:Elliptic}.
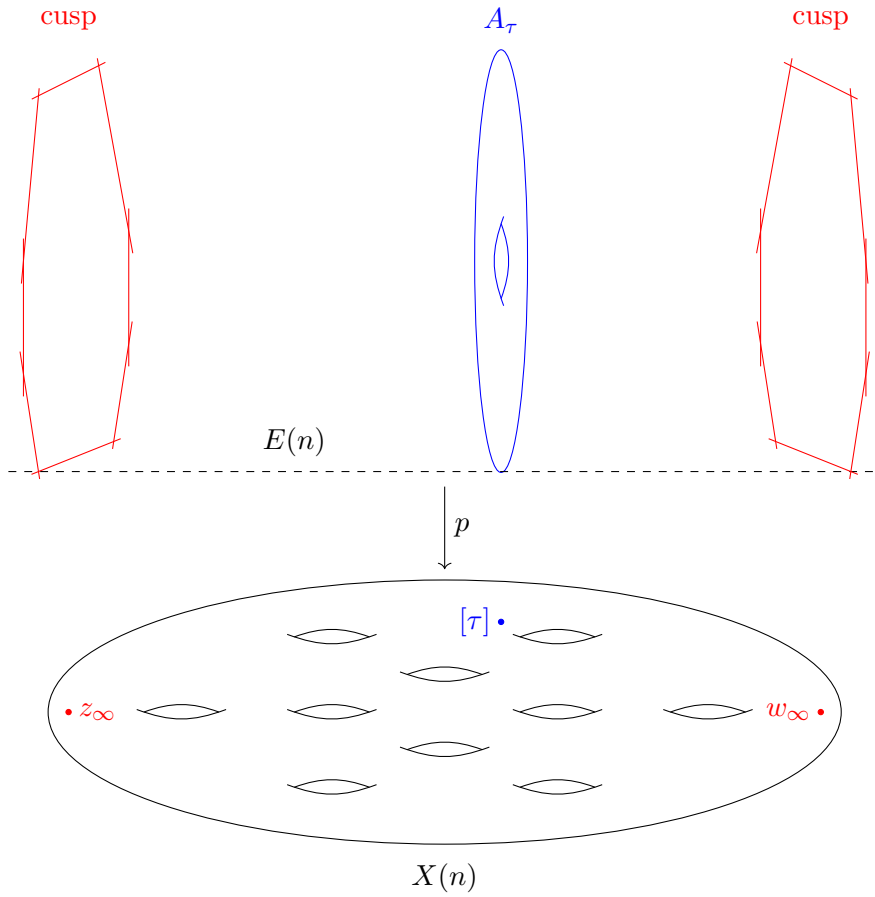
\begin{figure}[h]
\centering
\begin{tikzpicture}
\draw (0,0) ellipse (150pt and 50pt);
\node at (0,-2.2) {$X(n)$};
\draw [shorten >=-0.1cm,shorten <=-0.1cm] (1,1) to [bend right=20] (2,1);
\draw (1,1) to [bend left=20] (2,1);
\draw [shorten >=-0.1cm,shorten <=-0.1cm] (-2,1) to [bend right=20] (-1,1);
\draw (-2,1) to [bend left=20] (-1,1);
\draw [shorten >=-0.1cm,shorten <=-0.1cm] (1,-1) to [bend right=20] (2,-1);
\draw (1,-1) to [bend left=20] (2,-1);
\draw [shorten >=-0.1cm,shorten <=-0.1cm] (-2,-1) to [bend right=20] (-1,-1);
\draw (-2,-1) to [bend left=20] (-1,-1);
\draw [shorten >=-0.1cm,shorten <=-0.1cm] (-0.5,0.5) to [bend right=20] (0.5,0.5);
\draw (-0.5,0.5) to [bend left=20] (0.5,0.5);
\draw [shorten >=-0.1cm,shorten <=-0.1cm] (-0.5,-0.5) to [bend right=20] (0.5,-0.5);
\draw (-0.5,-0.5) to [bend left=20] (0.5,-0.5);
\draw [shorten >=-0.1cm,shorten <=-0.1cm] (-2,0) to [bend right=20] (-1,0);
\draw (-2,0) to [bend left=20] (-1,0);
\draw [shorten >=-0.1cm,shorten <=-0.1cm] (-4,0) to [bend right=20] (-3,0);
\draw (-4,0) to [bend left=20] (-3,0);
\draw [shorten >=-0.1cm,shorten <=-0.1cm] (3,0) to [bend right=20] (4,0);
\draw (3,0) to [bend left=20] (4,0);
\draw [shorten >=-0.1cm,shorten <=-0.1cm] (1,0) to [bend right=20] (2,0);
\draw (1,0) to [bend left=20] (2,0);
\draw [blue, fill=blue] (0.75,1.2) circle (1pt) node [left] {$[\tau]$};
\draw [blue] (0.75, 6) ellipse (10pt and 80pt);
\draw [blue, shorten >=-0.1cm,shorten <=-0.1cm] (0.75, 6.5) to [bend right=20] (0.75, 5.5);
\draw [blue] (0.75,6.5) to [bend left=20] (0.75,5.5);
\node [blue] at (0.75, 9.2) {$A_\tau$};
\node [red] at (-5, 9.2) {cusp};
\node [red] at (5, 9.2) {cusp};
\draw (0,3) [->] to node [right] {$p$} (0,1.9);
\node at (-2,3.6) {$E(n)$};
\draw [dashed] (-5.8,3.2) -- (5.8,3.2);
\draw [red, fill=red] (-5,0) circle (1pt) node [right] {$z_\infty$};
\draw [red, shorten >=-0.1cm,shorten <=-0.1cm] (-5.4,8.2) -- (-4.6,8.6);
\draw [red, shorten >=-0.3cm,shorten <=-0.1cm] (-5.4,8.2) -- (-5.6,6);
\draw [red, shorten >=-0.3cm,shorten <=-0.1cm] (-4.6,8.6) -- (-4.2,6.4);
\draw [red, shorten >=-0.3cm,shorten <=-0.3cm] (-5.6,6) -- (-5.6,4.5);
\draw [red, shorten >=-0.3cm,shorten <=-0.3cm] (-4.2,6.4) -- (-4.2,4.9);
\draw [red, shorten >=-0.1cm,shorten <=-0.3cm] (-5.6,4.5) -- (-5.4,3.2);
\draw [red, shorten >=-0.1cm,shorten <=-0.3cm] (-4.2,4.9) -- (-4.4,3.6);
\draw [red, shorten >=-0.1cm,shorten <=-0.1cm] (-5.4,3.2) -- (-4.4,3.6);
\draw [red, fill=red] (5,0) circle (1pt) node [left] {$w_\infty$};
\draw [red, shorten >=-0.1cm,shorten <=-0.1cm] (5.4,8.2) -- (4.6,8.6);
\draw [red, shorten >=-0.3cm,shorten <=-0.1cm] (5.4,8.2) -- (5.6,6);
\draw [red, shorten >=-0.3cm,shorten <=-0.1cm] (4.6,8.6) -- (4.2,6.4);
\draw [red, shorten >=-0.3cm,shorten <=-0.3cm] (5.6,6) -- (5.6,4.5);
\draw [red, shorten >=-0.3cm,shorten <=-0.3cm] (4.2,6.4) -- (4.2,4.9);
\draw [red, shorten >=-0.1cm,shorten <=-0.3cm] (5.6,4.5) -- (5.4,3.2);
\draw [red, shorten >=-0.1cm,shorten <=-0.3cm] (4.2,4.9) -- (4.4,3.6);
\draw [red, shorten >=-0.1cm,shorten <=-0.1cm] (5.4,3.2) -- (4.4,3.6);
\end{tikzpicture}
\caption{A compact elliptic fibration $E(n)$}\label{fig:Elliptic}
\end{figure}

Livn\'e showed that $E(9)$ can be realized as the underlying space for a nonarithmetic compact complex hyperbolic orbifold $\Gam \bs \bbB^2$. Since $X(9)$ has genus $10$, the projection $p : E(9) \to X(9)$ induces a surjective homomorphism from $\Gam$ onto a cocompact Fuchsian group. In particular, $\Gam$ has a nontrivial universal homomorphism. Taking $\Gam$ and its finite index subgroups proves the theorem.
\end{pf}

\begin{rem}\label{rem:DMexs}
Some nonarithmetic Deligne--Mostow lattices \cite{DeligneMostow} also admit homomorphisms onto cocompact Fuchsian groups from \emph{forgetful maps} associated with contracting weights. See for example \cite[Thm.\ 3.1(iv)]{Deraux}.
\end{rem}

\subsection{Automorphisms of Fuchsian groups}\label{ssec:FuchsianAut}

The purpose of this subsection is to prove \Cref{thm:Fixpsubs}.

\begin{pf}[Proof of \Cref{thm:Fixpsubs}]
Let $p$ be a prime, $T$ be a surface group of genus $g \ge2$, and $\phi \in \Aut(T)$ have the property that $\phi(N) = N$ for every finite index normal subgroup with index a power of $p$. Consider the pro-$p$ completion $T_{\wh{p}}$ of $T$. Then $\varphi$ extends to a continuous automorphism of $T_{\wh{p}}$ that is \emph{normal}, i.e., so that $\varphi(\wh{N})$ for every closed normal subgroup of $T_{\wh{p}}$. Indeed, every such normal subgroup is an intersection of finite index normal subgroups of $T_{\wh{p}}$, and finite index normal subgroups of $T_{\wh{p}}$ are in natural one-to-one correspondence with the normal $p$-power index subgroups of $T$.

Then $T_{\wh{p}}$ is topologically finitely presented, since $T$ is finitely presented, hence $T_{\wh{p}}$ is \emph{psuedo $p$-free} in the sense of Jarden--Ritter (see \cite[\S 1]{JardenRitter} for the definition) by \cite[Thm. 4]{JardenRitter}. Thus the extension of $\varphi$ to $T_{\wh{p}}$ is inner by \cite[Thm.\ 1]{JardenRitter}. Moreover, Paris proved that $T$ is $p$-conjugacy separable \cite[Thm.\ 1.7]{Paris}, which implies that $\varphi$ must preserve conjugacy classes of $T$ (i.e., it is \emph{pointwise inner}). Therefore $\varphi$ is a normal automorphism of $T$, meaning that it preserves all normal subgroups of $T$. Finally, it is a theorem of Bogopolski, Kudryavtseva, and Zieschang that normal automorphisms of higher-genus surface groups are inner \cite[Thm.\ 1.3]{BKZ}, which completes the proof of the theorem.
\end{pf}

\section{The main result}\label{sec:Main}

Suppose that $\Lam < \PU(2,1)$ is a cocompact nonarithmetic lattice with nontrivial universal homomorphism $\rho$. Following \Cref{ssec:PU21}, let $\Gam < \PU(2,1)$ be the unique maximal lattice containing $\Lam$. By \Cref{lem:TFFuchsian} and by passing to a finite index subgroup if necessary, we can assume that the universal homomorphism for $\Lam$ has a torsion-free factor and that $\Lam$ is a characteristic subgroup of $\Gam$. If $M_{\mathrm{tf}} \coloneqq \ker(\rho_{\mathrm{tf}})$ denotes the kernel of the torsion-free universal homomorphism for $\Lam$, then $M_{\mathrm{tf}}$ is an infinite index characteristic subgroup of $\Gam$ by \Cref{lem:UniversalChar} and because a characteristic subgroup of a characteristic subgroup is characteristic in the larger group. Set $T \coloneqq \Lam / M_{\mathrm{tf}}$, so $T = \rho_{\mathrm{tf}}(\Lam)$.

As described in \Cref{ssec:PU21}, away from finitely many prime ideals $\frakp$ of the ring $\Tr\Ad \Gam$, there is a well-defined congruence quotient $\wh{Q} \coloneqq \Gam / \Gam(\frakp)$ with trivial center and commutator subgroup $Q$ that is a nonabelian finite simple group. Moreover, $\Lam$ and $M_{\mathrm{tf}}$ are Zariski dense in $\PU(2,1)$, hence after possibly passing to a finite index subgroup of $\Lam$, we can assume that they both map onto $Q$ for infinitely many $\frakp$. Fix one such $\frakp$ for the remainder of this section. Then define $\Lam(\frakp) \coloneqq \Gam(\frakp) \cap \Lam$ and $M_{\mathrm{tf}}(\frakp) \coloneqq \Gam(\frakp) \cap M_{\mathrm{tf}}$, noting that both $\Lam(\frakp)$ and $M_{\mathrm{tf}}(\frakp)$ are normal subgroups of $\Gam$, since each is the intersection of normal subgroups of $\Gam$.

\begin{lem}\label{lem:QuoSetup}
With the notation established in this section,
\begin{align*}
\rho_{\mathrm{tf}}(\Lam(\frakp)) &= T \\
\Lam / M_{\mathrm{tf}}(\frakp) &\cong Q \times T \\
\Gam / M_{\mathrm{tf}}(\frakp) &\cong \wh{Q} \times \wh{T}
\end{align*}
where $\wh{T} \coloneqq \Gam / M_{\mathrm{tf}}$ contains $T$ as a finite index subgroup.
\end{lem}
\begin{pf}
The first equality holds since $M_{\mathrm{tf}}$ contains coset representatives for $\Lam(\frakp)$ in $\Lam$. The second holds since $M_{\mathrm{tf}} / M_{\mathrm{tf}}(\frakp)$ and $\Lam(\frakp) / M_{\mathrm{tf}}(\frakp)$ are disjoint normal subgroups of $\Lam / M_{\mathrm{tf}}(\frakp)$, and the third is similar.
\end{pf}

\begin{rem}
Note that $T$ is then a product of factors of the torsion-free universal homomorphism of $\Lam(\frakp)$, but there may be additional factors.
\end{rem}

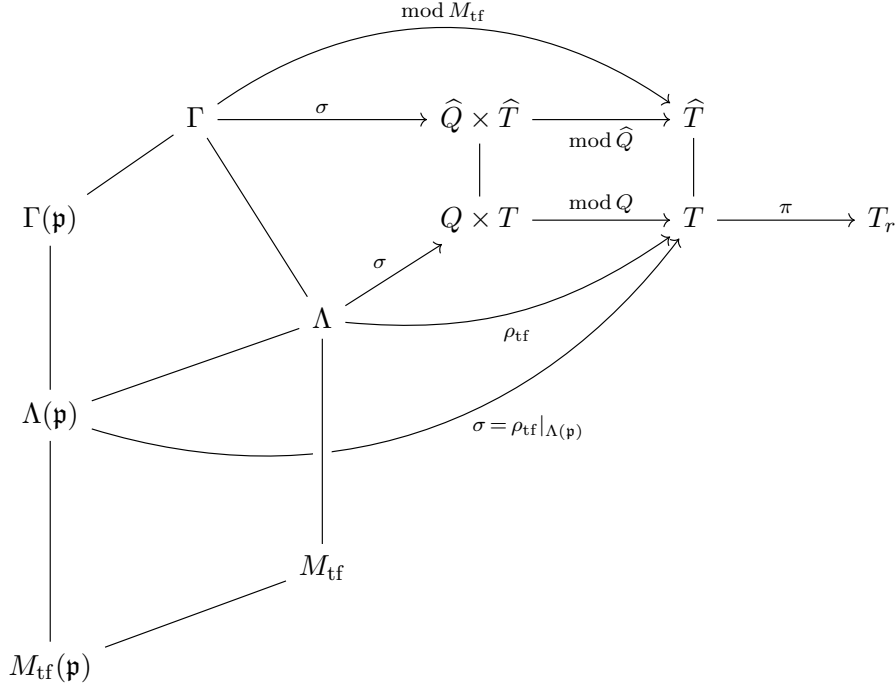
\begin{figure}[h]
\centering
\begin{tikzcd}
& \Gam \arrow[ddr, -] \arrow[dl, -] \arrow[rr, ->, "\sig"] \arrow[rrrr, bend left = 35, ->, "\textrm{mod}\, M_{\mathrm{tf}}"] & & \wh{Q} \times \wh{T} \arrow[d, -] \arrow[rr, ->, "\textrm{mod}\, \wh{Q}" below] & & \wh{T} \arrow[d, -] \\
\Gam(\frakp) \arrow[dd, -] & & & Q \times T \arrow[rr, ->, "\textrm{mod}\, Q"] && T \arrow[rr, ->, "\pi"] & & T_r \\
& & \Lam \arrow[dll, -] \arrow[ur, ->, "\sig"] \arrow[urrr, bend right=20, ->, "\rho_{\mathrm{tf}}" {xshift=0.75em, yshift=-1.25em}] \\
\Lam(\frakp) \arrow[ddd, -] \arrow[uurrrrr, ->, bend right=35, "\sig\, =\, \rho_{\mathrm{tf}}|_{\Lam(\frakp)}" {xshift=6em, yshift=-0.5em}] & & \\
& & \\
& & M_{\mathrm{tf}} \arrow[dll, -] \arrow[uuu, -, crossing over] & \\
M_{\mathrm{tf}}(\frakp) & &
\end{tikzcd}
\caption{The relationships between various groups}\label{fig:BigDiagram}
\end{figure}

In what follows, $\pi : T \to T_r$ will denote projection onto the $r^{th}$ surface group factor of $T$ for an arbitrary $r$ that will be fixed for the remainder of this section. Let $\sig$ denote projection of $\Gam$ onto $\Gam / M_{\mathrm{tf}}(\frakp)$ and also the restriction of $\sig$ to any finite index subgroup of $\Gam$. This leads to the diagram of groups and surjective homomorphisms given in \Cref{fig:BigDiagram}.
Define
\begin{equation}\label{eq:TauDef}
\tau \coloneqq \pi \circ \rho_{\mathrm{tf}} : \Lam \lra T_r
\end{equation}
and $M \coloneqq \ker(\tau|_\Del)$. For a finite index subgroup $\Del \le \Gam$, let $N_\Gam(\Del)$ denote its normalizer in $\Gam$. If $\Del$ contains $M$, then any $\gam \in N_\Gam(\Del)$ normalizing $M$ induces an automorphism of $\Del / M$ through $\pi \circ \rho_{\mathrm{tf}}$. The following proposition/construction is essential for proving \Cref{thm:Main}.

\begin{prop}\label{prop:FindD}
There is a finite index subgroup $\Del \le \Lam$ containing $M_{\mathrm{tf}}(\frakp)$ so that the group
\begin{equation}\label{eq:DDef}
D \coloneqq\!\left\{\gam \in N_\Gam(\Del)\ :\ \exists\ t_\gam \in T_r \textrm{ so } \tau(\gam \de \gam^{-1}) = t_\gam \tau (\de) t_\gam^{-1}\ \forall\ \de \in \Del \right\}
\end{equation}
consisting of the elements of $N_\Gam(\Del)$ that normalize $M$ and induce an inner automorphism of $T_r = \Del / M$ via $\tau$ is torsion-free.
\end{prop}
\begin{pf}
Since $Q$ is $2$-generated (see Window 2 of \cite{LubotzkySegal}) and $T_r$ is a higher-genus surface group, there is a surjective homomorphism $\psi_0 : T_r \to Q$. Define $\psi \coloneqq \psi_0 \circ \pi : T \to Q$ and the subgroup
\begin{equation}\label{eq:SDef}
S \coloneqq\!\left\{(\psi(t), t)\ :\ t \in T\right\}
\end{equation}
of $Q \times T$, which has finite index since $Q \times \{1\}$ forms a complete system of cosets for $S$ in $Q \times T$. Note that $\pi$ extends to $Q \times T$ by projecting away from $Q$, inducing a surjection from $S$ onto $T_r$, and define $\Del \coloneqq \sig^{-1}(S) \le \Lam$ with $D$ now defined as in the statement of the proposition.

Note that elements of $S$ are then all of the form
\begin{equation}\label{eq:Selts}
\sig(\de) = (\psi(\rho_{\mathrm{tf}}(\de)), \rho_{\mathrm{tf}}(\de)) \in Q \times T
\end{equation}
for an appropriate $\de \in \Del$. Given $\gam \in \Gam$, set $\sig(\gam) = (a_\gam, b_\gam) \in \wh{Q} \times \wh{T}$. If $\gam \in N_\Gam(\Del)$, then $\sig(\gam)$ normalizes $S$, thus
\begin{align}
a_\gam \psi(\rho_{\mathrm{tf}}(\de)) a_\gam^{-1} &= \psi(b_\gam \rho_{\mathrm{tf}}(\de) b_\gam^{-1}) \nonumber \\
&= \psi(\rho_{\mathrm{tf}}(\gam \de \gam^{-1})) \label{eq:1stConj}
\end{align}
for all $\de \in \Del$. If $\gam \in D$, then we additionally have
\begin{align}
\psi(\rho_{\mathrm{tf}}(\gam \de \gam^{-1})) &= \psi_0( \tau(\gam \de \gam^{-1})) \nonumber \\
&= \psi_0(t_\gam \tau(\de) t_\gam^{-1}) \label{eq:2ndConj}
\end{align}
for the appropriate $t_\gam \in T_r$. If $\gam$ were to be a nontrivial torsion element of $D$, then $t_\gam$ would act on $T_r$ with finite order. However, this implies that $t_\gam$ is trivial, since $T_r$ is torsion-free and has trivial center. Thus \Cref{eq:1stConj} and \Cref{eq:2ndConj} combined with triviality of $t_\gam$ imply that
\begin{equation}\label{eq:Dtf}
a_\gam \psi(\rho_{\mathrm{tf}}(\de)) a_\gam^{-1} = \psi(\rho_{\mathrm{tf}}(\de))
\end{equation}
for all $\de \in \Del$ and nontrivial torsion elements $\gam \in D$.

Since $\Gam(\frakp)$ is torsion-free, $a_\gam \in \wh{Q}$ is nontrivial for every torsion element $\gam \in \Gam$. Since $\wh{Q}$ has trivial center and the extension of $\psi$ to $Q \times T$ is surjective on $S$, there exists $\de \in \Del$ so that $\psi(\rho_{\mathrm{tf}}(\de))$ is nontrivial and does not commute with $a_\gam$. This implies that \Cref{eq:Dtf} cannot hold for all $\de \in \Del$. In particular, $D$ is torsion-free as claimed.
\end{pf}

For the remainder of this section, fix $\Del$ and $D$ from \Cref{prop:FindD}, and recall that $M$ is the kernel of $\tau|_\Del : \Del \to T_r$. Recall that $D$ consists of precisely the elements of $N_\Gam(\Del)$ that preserve $M$ and induce an inner automorphism of $T_r$ through $\tau|_\Del$. Equivalently, $D = C \Del$, where
\begin{equation}\label{eq:CDef}
C \coloneqq\!\left\{\gam \in N_\Gam(\Del)\ :\ \tau(\gam \de \gam^{-1}) = \tau (\de)\ \textrm{for all}\ \de \in \Del \right\}
\end{equation}
is the subgroup of $N_\Gam(\Del)$ inducing the trivial automorphism.

The remainder of this section follows arguments in \cite[\S 4]{BelolipetskyLubotzky} fairly closely, though some points remain delicate due to the fact that $\Del$ is not necessarily normal in $\Gam$. Varying the primes $p_j$ in the following proposition will be necessary for the `infinitely many' part of \Cref{thm:Main}.

\begin{prop}\label{prop:ChoosepSubgps}
Let $d+1 = [\Gam : D]$, choose primes $p_0 < \cdots < p_d$ all greater than $\max\!\left\{|G|, |\Aut(G)|\right\}^2$, and fix coset representatives $\gam_0, \dots, \gam_d$ for $D$ in $\Gam$ with $\gam_0 = 1$. For each $0 \le j \le d$ there are subgroups $K_j \triangleleft \Del$ with the following properties:
\begin{itemize}

\item[$\star$] $[\Del : K_j]$ is a nontrivial power $p_j^{\al_j}$ of $p_j$ for all $0 \le j \le d$;

\item[$\star$] $M \le K_j$ for all $0 \le j \le d$;

\item[$\star$] $\gam_j K_j \gam_j^{-1} \neq K_j$ for $1 \le j \le d$;

\item[$\star$] further, $\wt{K}_j \coloneqq \langle K_j, M_\mathrm{tf} \rangle \le \Lam$, and $\gam_j \wt{K}_j \gam_j^{-1} \neq \wt{K}_j$ for $j \ge 1$.

\end{itemize}
Moreover, $\wt{K}_j$ is equivalently the preimage in $\Lam$ of $\tau|_\Del(K_j) \le T_r$.
\end{prop}
\begin{pf}
If $d = 0$, then the preimage in $\Del$ of any normal subgroup of $T_r$ of index $p_0$ suffices. If $d > 0$ and there were no such $K_j$ for some $j > 0$, then $\gam_j$ would normalize every $p_j$-power normal subgroup of $\Del$ containing $M$. Since surface groups are residually $p_j$ (indeed, they are even residually free \cite[Thm.\ 1]{Baumslag}), the intersection of all normal $p_j$-power index subgroups of $\Del$ containing $M$ equals $M$, and thus $\gam_j$ normalizes $M$ and induces an automorphism of $T_r$ through $\tau|_\Del$. Since $\gam_j \notin D$, this induced automorphism is not inner. Therefore, \Cref{thm:Fixpsubs} implies that there is some $K_j$ not normalized by $\gam_j$.

Now, since $M_\mathrm{tf}$ is normal in $\Gam$, we have that
\begin{equation}\label{eq:Mquo}
M_\mathrm{tf} / M_\mathrm{tf}(\frakp) = Q \times \{1\} \le Q \times T
\end{equation}
contains coset representatives for $S$ in $Q \times T$, which moreover implies that $M_\Lam \coloneqq \langle M, M_\mathrm{tf}\rangle$ is the kernel of $\tau$. One equivalently checks that
\begin{equation}\label{eq:LamM}
\tau(M_\Lam) = Q \times \{1\} \times \prod_{\ell \neq r} T_\ell
\end{equation}
is the image in $Q \times T$ of the kernel of $\tau$. Then $\wt{K}_j$ is the preimage in $\Lam$ of $\tau|_\Del(K_j) < T_r$, giving the final claim of the proposition. Moreover, $\gam_j$ normalizes $M_\Lam$ for the exact same reason as for $M$, so we conclude that $\wt{K}_j$ is also not normalized by $\gam_j$. This concludes the proof.
\end{pf}

In fact, the proof of \Cref{prop:ChoosepSubgps} produces many choices of $K_j$. Fixing $p_0, \dots, p_d$ and $K_j \triangleleft \Del$ of index $p_j^{\al_j}$ from \Cref{prop:ChoosepSubgps}, define
\begin{align}
K &\coloneqq \bigcap_{j = 0}^d K_j \label{eq:KDef} \\
k &\coloneqq [\Del : K] \label{eq:kDef}
\end{align}
for the remainder of this section. If $T_r$ is a surface group of genus $h$, then the image of $K$ in $T_r$ under $\tau$ is a surface group of genus $h_K \coloneqq 1 + k(h - 1)$. Note that our hypotheses on the primes $p_j$ ensure that $h_K > 2 \log_2\!|G|$. A preliminary result, whose proof simply bootstraps the analogous calculation for the free group in \cite[\S 4]{BelolipetskyLubotzky}, is required.

Given a nontrivial finite group $G$, let $\calM_G$ be the collection of all normal subgroups $A \triangleleft K$ containing $M$ with $K/A \cong G$. We first prove that $\calM_G$ is quite large.

\begin{lem}\label{lem:SurfaceToG}
For any nontrivial finite group $G$, $\calM_G$ contains at least $|G|^{h_K - 2 \log_2 |G|}$ elements.
\end{lem}
\begin{pf}
Since $G$ can be generated by $\log_2\! |G|$ elements \cite[Lem.\ 1.2.2]{LubotzkySegal}, and $K/M \le T_r$ has a free quotient $F$ of rank $h_K$, which is greater than $2 \log_2\! |G|$ by our assumption on the primes $p_j$, one can send the first $\log_2\! |G|$ generators of $F$ to a fixed generating set for $G$ then do anything with the remaining $h_k - \log_2\! |G|$ generators. Two homomorphisms $\varphi_1, \varphi_2 : K \to G$ have the same kernel if and only if $\varphi_2 = \beta \circ \varphi_1$ for some $\beta \in \Aut(G)$. Also, ${|\Aut(G)| \le |G|^{\log_2\! |G|}}$ since any automorphism of $G$ is determined by what it does to the generating set of size at most $\log_2\! |G|$. Therefore, there are at least
\begin{equation}\label{eq:AutSizeBound}
\frac{1}{|\Aut(G)|} |G|^{h_K - \log_2\! |G|} \ge |G|^{h_K - 2 \log_2\! |G|}
\end{equation}
distinct normal subgroups $A \triangleleft K$ with $K/A \cong G$.
\end{pf}

We now show that elements of $\calM_G$ have normalizer contained in $D$.

\begin{prop}\label{prop:NormInD}
If $G$ is nontrivial and $A \in \calM_G$, then $N_\Gam(A) \le D$.
\end{prop}
\begin{pf}
Recall that $D = \Gam$ if $d = 0$, so there is nothing to prove in that case. Otherwise, suppose that $\gam \in N_\Gam(A) \ssm D$. Notice that $A$ is precisely the preimage in $\Del$ of $\conj{A} \coloneqq \tau(A) < T_r$, since $A$ contains $M$ by hypothesis. Then, the preimage $\wt{A}$ of $\conj{A}$ in $\Lam$ under projection of $\Lam$ onto $T_r$ equals $\langle A, M_\mathrm{tf} \rangle$. Indeed, $\wt{A}$ is equivalently the preimage in $\Lam$ of
\begin{equation}\label{eq:wtADef}
Q \times \conj{A} \times \prod_{\ell \neq r} T_\ell \le Q \times T
\end{equation}
under reduction modulo $M_\mathrm{tf}(\frakp)$, and $M_\mathrm{tf} / M_\mathrm{tf}(\frakp) = Q \times \{1\} \le Q \times T$ then contains the elements of $Q \times T$ that along with
\begin{equation}\label{eq:ImageOfA}
A / M_\mathrm{tf}(\frakp) =\! \left\{(\psi(t), t)\ :\ t \in \conj{A} \times \prod_{\ell \neq r} T_\ell \right\}\! \le Q \times T
\end{equation}
generate the left-hand side of \Cref{eq:wtADef}.

Since $\gam$ normalizes $A$ and $M_\mathrm{tf}$, it also normalizes $\wt{A}$. From \Cref{prop:ChoosepSubgps}, we then similarly have $\wt{K}_j \coloneqq \langle K_j, M_\mathrm{tf} \rangle \triangleleft \Lam$, which is the preimage in $\Lam$ of $\conj{K}_j \coloneqq \tau(K_j)$, and
\begin{equation}\label{eq:wtKDef}
\wt{K} \coloneqq \bigcap_{\ell = 0}^d \wt{K}_\ell
\end{equation}
equals $\langle K, M_\mathrm{tf} \rangle$, since it is also the preimage of $\bigcap \conj{K}_j \le T_r$ in $\Lam$. Since $\wt{K}_j$ is the unique maximal $p_j$-power index subgroup of $\Lam$ containing $\wt{A}$, it follows that $\gam$ also normalizes $\wt{K}_j$ for all $j$.

Write $\gam = \gam_i \de$ for some $\de \in D$ and $\gam_i$ one of the coset representatives for $D$ in $\Gam$ fixed for this section. Then $\de$ normalizes $M$, $M_\mathrm{tf}$, and $\Lam$, so it also normalizes the kernel $M_\Lam$ of $\tau$, since $M_\Lam = \langle M, M_\mathrm{tf}\rangle$ as in the proof of \Cref{prop:ChoosepSubgps}. Since $\de$ acts on $T_r$ by an inner automorphism, it also normalizes $\wt{K}_j$. It follows that $\gam_i = \gam \de^{-1}$ normalizes each $\wt{K}_j$, hence the case $j = i$ contradicts \Cref{prop:ChoosepSubgps} and thus proves the proposition.
\end{pf}

The next proposition is the penultimate technical result needed in the case where $G$ is nontrivial.

\begin{prop}\label{prop:FindA}
Suppose $G$ is a nontrivial finite group and $\calM_G$ is defined as in \Cref{prop:NormInD} with respect to $K \triangleleft \Del$ as defined in \Cref{eq:KDef}. After possibly increasing the smallest prime $p_0$ used to define $K$, there exists at least one $A$ in $\calM_G$ so that:
\begin{align}
N_D(A) &= N_\Gam(A) \label{eq:ND=NGam} \\
N_\Del(A) / A &\cong G \label{eq:NDel(A)/A=G} \\
N_D(A) / A &\cong G \times C / M \label{eq:ND(A)/A}
\end{align}
\end{prop}
\begin{pf}
Note that $N_D(A) = N_\Gam(A)$ for all $A \in \calM_G$ by \Cref{prop:NormInD}, so \Cref{eq:ND=NGam} always holds. Now, observe that $KC \le N_\Gam(A)$. To see this, first note that $K$ normalizes $A$ by definition. In addition, $C$ normalizes $\Del$ and $M$ by definition, and since its action on $T_r$ by conjugation is trivial, it also normalizes the image $\conj{A}$ of $A$ in $T_r$. Therefore, $C$ normalizes $A$. Note that $K$ and $C$ are both normal in $D$, so $KC$ is a group.

Recall that $D = \Del C$ with $C \cap \Del = M$, since $\Del / M = T_r$ has trivial center. Moreover, recall that $M \triangleleft D$, hence:
\begin{align}
D / M &= \Del / M \times C / M \label{eq:D/M} \\
N_D(A) / M &= N_\Del(A) / M \times C / M \label{eq:ND(A)/M}
\end{align}
Projection onto $\Del / M$ therefore implies that it suffices to prove that there is an $A \in \calM_G$ with $N_\Del(A) = K$. Indeed, this immediately gives \Cref{eq:NDel(A)/A=G} by definition of $\calM_G$ and \Cref{eq:ND(A)/A} then follows from \Cref{eq:ND(A)/M}. The existence of such an $A \in \calM_G$ follows from the same counting argument found on \cite[p.\ 466-467]{BelolipetskyLubotzky}. The proof is recalled for the reader's convenience.

The subgroup $N_\Del(A)/M$ of the genus $h$ surface group $T_r = \Del / M$ is a surface group of genus $h_N \coloneqq 1 + [\Del : N_\Del(A)](h-1)$. Then, recalling $k$ from \Cref{eq:kDef},
\begin{equation}\label{eq:AIndex}
n \coloneqq [N_\Del(A) / M : A / M] = [N_\Del(A) : A] = |G| \frac{k}{[\Del : N_\Del(A)]}
\end{equation}
since $M \le A$, and \cite[Thm.\ 2.7]{LubotzkySegal} then implies that $N_\Del(A) / M$ has at most
\begin{equation}\label{eq:NormalBound1}
n^{2 h_N \! \log_2(n) c} \le 2^{2 c h [\Del : N_\Del(A)]\! \log_2(n)^2}
\end{equation}
normal subgroups of index $n$, where $c$ is an absolute constant. Indeed, $N_\Del(A)/M$ can be generated by $2 h_N$ elements and the number of normal subgroups of a given index is bounded from above by the analogous count for the free group on the same number of generators.

Fix any $x$ so that $\max\{|G|, |\Aut(G)|\} < x < \sqrt{p_0}$. Since $K \le N_\Del(A)$ by definition of $A$, proper containment of $K$ in $N_\Del(A)$ along with the fact that $p_0$ is the smallest divisor of $k$ by construction of $K$, implies that
\begin{equation}\label{eq:IndexBound}
[\Del : N_\Del(A)] \le \frac{k}{p_0} \le \frac{k}{x^2} < \frac{k}{x}
\end{equation}
with $k \ge x^2$ by construction. Incorporating this into \Cref{eq:NormalBound1}, a given finite index subgroup $L \le \Del$ containing $M$ and properly containing $K$ can normalize at most
\begin{equation}\label{eq:NormalBound2}
2^{2 c h \frac{k}{x}\! \log_2(|G| x)^2}
\end{equation}
elements of $\calM_G$.

However, any such $L$ has the property that $L / K$ is a nontrivial subgroup of the finitely generated (in fact, finite) nilpotent group $\Del / K$, hence it is subnormal with index at most $k / x$. Thus \cite[Thm.\ 2.3]{LubotzkySegal} implies that there are at most $2^{2 h k/x}$ possible $L$. Thus there are at least
\begin{equation}\label{eq:NormalBound3}
|G|^{k(h - 1) - 2 \log_2(|G|)}
\end{equation}
possibilities for $A$ and at most
\begin{equation}\label{eq:NormalBound4}
2^{(h k / x)(1 + c \log_2(|G| x)^2)}
\end{equation}
possibilities have normalizer strictly larger than $K$. The inequality
\begin{equation}\label{eq:NormalBound5}
\frac{2 c}{x} \log_2(|G| x)^2 + \frac{(\log_2 |G|)^2}{k} + \frac{2}{x} < \log_2 |G|
\end{equation}
then implies that the quantity in \Cref{eq:NormalBound3} is greater than the quantity in \Cref{eq:NormalBound4}. Making $p_0$ sufficiently large forces $k$ to be large, so we can assume $(\log_2 |G|)^2 / k < (\log_2 |G|) / 3$, and allows us to choose $x$ so that $2 c \log_2(|G| x)^2 / x$ and $2 / x$ are also bounded from above by $(\log_2 |G|) / 3$. In particular, if $p_0$ is sufficiently large, then there is an $x$ satisfying the above constraints that is sufficiently large that \Cref{eq:NormalBound5} holds. Therefore there are $A \in \calM_G$ so that $N_\Del(A) = K$, as desired. This completes the proof of the proposition.
\end{pf}

The following is the denouement to the sequence of results needed to prove \Cref{thm:Main} in the case where $G$ is a nontrivial finite group.

\begin{prop}\label{prop:FinalNontrivialG}
Suppose $G$ is a nontrivial finite group, $A \in \calM_G$ satisfies the conclusions of \Cref{prop:FindA}, and $B = AC$ with $C$ as in \Cref{eq:CDef}. Then $N_\Gam(B) = N_\Gam(A)$. In particular, $N_\Gam(B) / B \cong G$.
\end{prop}
\begin{pf}
Since $C$ is normal in $D$, it follows that $N_\Gam(A) \le N_\Gam(B)$. Indeed,
\begin{equation}\label{eq:NDA}
N_D(A) \le N_D(AC) = N_D(B) \le N_\Gam(B)
\end{equation}
by \Cref{prop:NormInD}. The opposite inclusion then proves the proposition, as the second conclusion subsequently follows from \Cref{eq:ND=NGam} and \Cref{eq:ND(A)/A}. Moreover, since $\Del \cap C = M$ by \Cref{eq:D/M} and $M \le A$, it follows that $\Del \cap B = A$. In particular, it suffices to show that $\gam A \gam^{-1} \le \Del$ for all $\gam \in N_\Gam(B)$.

Fix $g \in N_\Gam(B)$, and let $\sig(g) = (g_1, g_2) \in \wh{Q} \times \wh{T}$ denote its reduction modulo $M_\mathrm{tf}(\frakp)$. Recall that the image of $\Del$ in $\wh{Q} \times \wh{T}$ is the group $S$ defined in \Cref{eq:SDef}. Since $C \le N_\Gam(\Del)$ by definition, the image $(a_\gam,b_\gam) \in \wh{Q} \times \wh{T}$ of an element $\gam \in C$ must satisfy \Cref{eq:1stConj} and \Cref{eq:2ndConj} with $t_\gam = 1$. Then $\psi$ factors through the projection $\pi$ onto $T_r$, so $a_\gam = 1$, since $\wh{Q}$ has trivial center. Moreover, the element $b_\gam \in \wh{T}$ must be in the centralizer $Z_{\wh{T}}(T_r)$ of the direct factor $T_r$ of $T$ in order to act trivially through $\pi$.

Now, $ac \in B$ satisfies $\sig(ac) = (\psi(t_a), t_a h_c)$ for some $t_a \in \rho_{\mathrm{tf}}(A)$ and $h_c \in Z_{\wh{T}}(T_r)$. If $a \in A$, then $g a g^{-1} = a^\prime c \in AC$, which gives
\begin{align}
\sig\!\left(g a g^{-1}\right) &= \left(\psi(t^\prime), t^\prime h\right) \nonumber \\
&= \left(g_1 \psi(t) g_1^{-1}, g_2 t g_2^{-1}\right) \label{eq:ACendgame}
\end{align}
for the appropriate $t, t^\prime \in T$ and $h \in Z_{\wh{T}}(T_r)$. Then $g_2 \in \wh{T}$ normalizes $T$, so $t^\prime h \in T$ and then $h \in T$. However, this implies that
\begin{equation}\label{eq:hCentralizes}
h \in \{1\} \times \prod_{\ell \neq r} T_\ell
\end{equation}
in order to centralize the $T_r$ factor, hence $\psi(t^\prime h) = \psi(t^\prime)$ since $\psi$ factors through projection onto $T_r$, meaning precisely that $\sig(g a g^{-1}) \in S$. This proves that $g a g^{-1} \in \Del$ for all $g \in N_\Gam(B)$ and $a \in A$, which is what we needed to prove.
\end{pf}

\Cref{prop:FinalNontrivialG} completes the proof of \Cref{thm:Main} for nontrivial $G$. Indeed, we produced a finite index torsion-free subgroup $B \le \Gam$ such that $N_\Gam(B) / B \cong G$. The formal proof of \Cref{thm:Main} that closes this section describes more precisely how this conclusion gives the desired statement.

We now turn to the case of the trivial group. Fix $x$ and primes $\{p_j\}$ exactly as above, so
\begin{equation}\label{eq:Redefine}
1 \ll x < p_0 < \cdots < p_d
\end{equation}
where $d + 1 = [\Gam : D]$, $K_j < \Del$ are the subgroups provided by \Cref{prop:ChoosepSubgps}, and $K = \bigcap K_j$. The following result, which closely follows the analogue in \cite{BelolipetskyLubotzky}, is our version of \Cref{prop:FinalNontrivialG} for the case where $G$ is trivial.

\begin{prop}\label{prop:TrivialCase}
If $q$ is a sufficiently large prime, then there is an element $A$ in the set $\calM_q$ of subgroups $A \le K$ with index $q$ containing $M$ such that $N_\Gam(A) = AC$, and thus $B = AC$ is torsion-free with $N_\Gam(B) = B$.
\end{prop}
\begin{pf}
Let $q$ be a prime number greater than $p_0$ with order roughly $k$ that is not among the $\{p_j\}$ and $\calM_q$ be the all the subgroups of $K$ with index $q$ that contain $M$. If $q$ is sufficiently large, then
\begin{equation}\label{eq:SubgroupLowerBound}
|\calM_q| > (q!)^{2 h_K - 2}
\end{equation}
by \cite[Thm.\ 14.4.5]{LubotzkySegal}, recalling that $h_K$ is the genus of the surface group $\conj{K}$.

If $N_\Del(A)$ strictly contains $A$, then $[\Del : N_\Del(A)]$ strictly divides $k q$, hence
\begin{equation}\label{eq:TrivialBound1}
[\Del : N_\Del(A)] \le \frac{k q}{x}
\end{equation}
with $x$ from \Cref{eq:Redefine}. Let $L$ be any subgroup of $\Del$ containing $M$ with index properly dividing $k q$. Then $L / M$ is a surface group of genus
\begin{equation}\label{eq:Trivialhbound}
1 + [\Del : L](h - 1) \le [\Del : L] h
\end{equation}
and therefore \cite[Thm.\ 2.7]{LubotzkySegal} implies that there are less than
\begin{equation}\label{eq:TrivialBound2}
2^{2 c h [\Del : L] \log_2(k q / [\Del : L])^2}
\end{equation}
groups in $\calM_q$ normalized by $L$, where $c$ is the same absolute constant appearing in \Cref{eq:NormalBound1}. Choosing $x$ sufficiently large by selecting sufficiently large initial primes $\{p_j\}$, the maximum in \Cref{eq:TrivialBound2} is attained when $L$ has index $k q / x$, giving the upper bound $2^{2 c h (kq / x) \log_2(x)^2}$ for the number of subgroups in $\calM_q$ normalized by the given $L$.

Since $L / M$ is a subgroup of the surface group $\Del / M$ of genus $h$, there are at most
\begin{equation}\label{eq:TrivialBound3}
\left(\left(\frac{k q}{x}\right)!\right)^{2 h} < \left(\frac{k q}{x}\right)^{2 c^\prime h k q / x}
\end{equation}
such subgroups by \cite[Cor.\ 1.1.2]{LubotzkySegal}, where $c^\prime$ is another absolute constant. It follows that there are at most
\begin{equation}\label{eq:TrivialBound4}
2^{2 h (k q /x)(c^\prime \log_2(k q / x) + c \log_2(x)^2)}
\end{equation}
subgroups in $\calM_q$ whose normalizer in $\Del$ is strictly larger than $A$. Since $q > x$ and $q$ is close to $k$, it is not hard to show using Stirling's formula that the quantity in \Cref{eq:SubgroupLowerBound} is larger than \Cref{eq:TrivialBound4} for sufficiently large $x$, meaning that there must be some $A \in \calM$ with $N_\Del(A) = A$.

The proof of \Cref{prop:NormInD} does not use anywhere that $A \triangleleft K$, so we again see that $N_\Gam(A) \le D$. Set $B = AC$. It then follows exactly the same as in the nontrivial case that $N_\Gam(B) = B$. Since $B \le D$, which is torsion-free by \Cref{prop:FindD}, this completes the proof of the proposition.
\end{pf}

We are now prepared to prove \Cref{thm:Main}.

\begin{pf}[Proof of \Cref{thm:Main}]
The only input to begin this section of the paper was a cocompact nonarithmetic lattice $\Lam < \PU(2,1)$ with nontrivial universal homomorphism. Such examples exist by \Cref{thm:Examples}. If $G$ is a nontrivial group (resp.\ the trivial group), then \Cref{prop:FinalNontrivialG} (resp.\ \Cref{prop:TrivialCase}) produces a torsion-free lattice $B \le \Lam$ whose normalizer $N(B) < \PU(2,1)$, which equals $N_\Gam(B)$ for $\Gam$ the unique maximal lattice in $\PU(2,1)$ containing $B$ by the Margulis dichotomy, has $N(B) / B \cong G$. Mostow rigidity now implies that $G$ is the automorphism group (equivalently, holomorphic isometry group) of the smooth compact complex hyperbolic $2$-manifold $B \bs \bbB^n$. Varying the primes $\{p_j\}$ in our construction gives infinitely many lattices with covolume tending to infinity for any fixed finite group $G$, hence another application of Mostow rigidity implies that there are infinitely many isomorphism classes of lattices. This proves the theorem.
\end{pf}

\section*{Appendix: Automorphism groups of surfaces}

The purpose of this appendix is to give a short account regarding smooth projective surfaces with finite automorphism group. We break up the discussion according to the Enriques--Kodaira classification of surfaces by their Kodaira dimension. To begin, we note that a projective surface is \emph{minimal} if it does not contain a rational curve (i.e., $\bbP^1$) of self-intersection $-1$, \emph{birationality} is the equivalence relation defined by having isomorphic dense open subvarieties, and every smooth projective surface is birational to a minimal surface. For a basic reference on the geometry of surfaces and the Kodaira--Enriques classification, see \cite{Beauville, BPV}.

\subsubsection*{Kodaira dimension $-\infty$}

A surface is rational if it is birational to the complex projective plane $\bbP^2$, and the group of birational automorphisms of $\bbP^2$ is the \emph{Cremona group} $\mathrm{Cr}_2(\bbC)$. Any finite group of automorphisms of a rational surface then gives a finite subgroup of $\mathrm{Cr}_2(\bbC)$. As with the automorphism group $\PGL_3(\bbC)$ of $\bbP^2$, $\mathrm{Cr}_2(\bbC)$ has the \emph{Jordan property}: Every finite subgroup has an abelian subgroup of universally bounded index. See \cite[\S4, \S5, \S10]{DolgachevIskovskikh} for a complete classification. Thus there are many finite groups that cannot even be subgroups of the automorphism group of a rational surface.

The surfaces of Kodaira dimension $-\infty$ that are not rational are the \emph{ruled surfaces}, which are $\bbP^1$ bundles $\pi : X \to C$ over a curve $C$ of genus $g \ge 1$. If $X$ is not rational, it has automorphism group fitting into an exact sequence
\begin{equation}\label{eq:RuledAut}
\begin{tikzcd} 1 \arrow[r] & \Aut_C(X) \arrow[r] & \Aut(X) \arrow[r, "\pi^*"] & \Aut(C) \end{tikzcd}
\end{equation}
with $\Aut_C(X)$ the subgroup commuting with projection \cite[Lem.\ 6]{Maruyama}. Moreover, there are examples with $\Aut_C(X)$ naturally trivial; see \cite[\S 2]{Maruyama}. In particular, we believe that there should be irrational ruled surfaces with any finite automorphism group, where the isomorphism $\Aut(X) \cong G$ is simply induced by Greenberg's theorem and an isomorphism $\Aut(X) \cong \Aut(C)$. However, we could not find a precise statement to this effect in the literature.

\medskip

As mentioned above for rational surfaces, if any two surfaces $S$ and $S^\prime$ are birational, then $\Aut(S^\prime)$ is naturally a subgroup of the group of birational automorphisms of $S$. For the remainder of this appendix, we use the following result in order to assume that surfaces are minimal.

\begin{thmnonum}[Thm.\ V.19 \cite{Beauville}]\label{thm:ToMinimal}
The group of birational automorphisms of a nonruled minimal surface that is not rational coincides with its automorphism group.
\end{thmnonum}

\subsubsection*{Kodaira dimension $0$}

An abelian surface has infinite automorphism group, since it acts on itself by translation. Bielliptic surfaces are of the form $(E_1 \times E_2) / G$, where $G$ acts on $E_1$ with fixed points and on $E_2$ by translations. Thus any infinite order translation of $E_2$ commutes with $G$ and descends to an infinite order automorphism of the bielliptic surface. Thus the only surfaces of Kodaira dimension zero that can have finite automorphism group are K3 and Enriques surfaces.

Kond\={o} showed that an algebraic K3 surface $X$ with finite automorphism group can only have $\Aut(X)$ trivial, $\bbZ / 2$, $(\bbZ / 2)^2$, $S_3 \times \bbZ / 2$, or a extension of one of these groups by $\bbZ / n$ for some $n \le 66$ \cite[p.\ 2]{KondoK3}. Kond\={o} also showed that the only groups appearing for an Enriques surface with finite automorphism group are the dihedral group $D_4$ of order $8$, $D_4 \ltimes (\bbZ / 2)^4$, $S_4$, $S_5$, $S_4 \times \bbZ / 2$, or $N \ltimes (\bbZ / 2)^4$ with $N$ the unique nontrivial semidirect product $\bbZ / 4 \ltimes \bbZ / 5$ \cite[Main Thm.]{KondoEnriques}. As with rational surfaces, we see that the finite subgroups of $\Aut(X)$ are very much restricted in Kodaira dimension zero.

\subsubsection*{Kodaira dimension $1$}

Surfaces of Kodaira dimension $1$ are the \emph{properly elliptic fibrations}: surfaces fibered over a curve with general fiber an elliptic curve not of Kodaira dimension $0$ or $-\infty$. This includes the direct products $E \times C$ with $E$ an elliptic curve and $C$ genus at least two, which clearly have infinite automorphism group. There is again an exact sequence as in \Cref{eq:RuledAut}, and it may be possible to realize any finite group as the automorphism group of a minimal properly elliptic fibration. Again, we could not find a proof of this fact in the literature.

\subsubsection*{Surfaces of general type}

Complex hyperbolic manifolds are of general type \cite[\S V.20]{BPV}, hence \Cref{thm:Main} immediately implies \Cref{thm:ProjectiveVersion}. Surfaces of general type all have finite automorphism group, which is directly analogous to the classical Hurwitz bound $|\Aut(C)| \le 84(g-1)$ for a smooth projective curve of genus $g$. In fact, $|\Aut(X)|$ can be bounded from above in terms of the self-intersection of the first Chern class of the canonical bundle; see for example \cite{Xiao}. Thus it is natural that the volume (or Euler characteristic, by Gauss--Bonnet) goes to infinity with $|G|$ for our examples.

\begin{rem}\label{rem:Hurwitz}
The second author conjectured for compact complex hyperbolic $2$-manifolds that $|\Aut(X)| \le 288 e(X)$, where $e(X)$ is the Euler characteristic, and proved the conjecture when $X$ is arithmetic \cite{StoverHurwitz}. Moreover, equality should come from normal subgroups of a particular arithmetic lattice in $\PU(2,1)$. This is completely analogous to the Hurwitz bound for curves, where equality holds if and only if the lattice in $\PU(1,1)$ is a normal subgroup of the $(2,3,7)$ triangle group, which is arithmetic.
\end{rem}

\bibliography{ComplexFiniteAut}

\end{document}